\documentclass[12pt,reqno]{amsart}

\usepackage[arrow,matrix,curve]{xy}

\usepackage[dvips]{graphicx} 

\usepackage{amssymb, latexsym, amsmath, amscd, array, hyperref, epigraph
}

\usepackage{hyperref}

\theoremstyle{definition}

\numberwithin{equation}{section}

\newcommand\N {{\mathbb N}} 

\newcommand\R {{\mathbb R}}
\newcommand\Q {{\mathbb Q}}

\newcommand\Los{{\L}o{\'s}}

\newcommand\astr{{}^\ast\R}

\newcommand\astf {{}^{\ast}\hspace{-2.6pt}f}

\author{Mikhail G. Katz}\address{M. Katz, Department of Mathematics,
Bar Ilan University, Ramat Gan 52900 Israel}
\email{katzmik@macs.biu.ac.il}

\author{Luie Polev}\address{L. Polev, Department of Mathematics, Bar
Ilan University, Ramat Gan 52900 Israel}\email{luiepolev@gmail.com}

\begin{document}

\thispagestyle{empty}


\title [From Pythagoreans and Weierstrassians to TIC] {From
Pythagoreans and Weierstrassians to true infinitesimal calculus}

\begin{abstract}
In teaching infinitesimal calculus we sought to present basic concepts
like continuity and convergence by comparing and contrasting various
definitions, rather than presenting ``the definition" to the students
as a monolithic absolute.  We hope that this could be useful to other
instructors wishing to follow this method of instruction.  A poll run
at the conclusion of the course indicates that students tend to favor
infinitesimal definitions over~$\epsilon,\delta$ ones.

\end{abstract}

\maketitle

\tableofcontents

\section{Introduction}
\label{one}

\epigraph{L\"ubsen defined the differential quotient first by means of
the limit notion; but along side of this he placed (after the second
edition) what he considered to be the \emph{true infinitesimal
calculus} -- a mystical scheme of operating with infinitely small
quantities.  \cite[p.\;217]{Kl08}}

Starting from the assumption that multiple approaches to the same
concept can facilitate student learning, during the 2014-2015 academic
year the authors taught \emph{True Infinitesimal Calculus} (TIC) based
on \cite{Ke86} to about 120 freshmen.  A similar course was taught
during the 2015-2016 year to about 130 students, and is being taught
during the 2016-2017 year to 150 students.  Most of the students had
already seen the basic techniques of the calculus in their highschool
courses.

Keisler's book was reviewed in \cite{Bi77} but Bishop's position is
dictated by a broad opposition to all of classical mathematics as
spelled out four years earlier in his ``Schizophrenia'' text
\cite{Bi73}; see \cite{KK11d} for details.  A sympathetic historical
account of infinitesimals can be found in \cite{Al14}.

In an effort to quantify student attitudes toward (1) the
\emph{Epsilontik} and (2) infinitesimals, we ran a poll at the end of
the course.  The goal was to compare student reactions to the two
approaches, as well as to gauge the helpfulness of each approach in
their eyes.  A total of 84 students participated in the poll.

Two-thirds of the respondents to the poll felt that infinitesimal
definitions of three key calculus concepts helped them understand the
concept, while only one in seven felt that~$\epsilon,\delta$
definitions helped them understand the concept.

We refer to the approach used in the course as TIC to distinguish it
from the traditional \emph{Epsilon Delta Calculus} (EDC).  While EDC
is often referred to as infinitesimal calculus, the use of the
adjective \emph{infinitesimal} in that term is something of a dead
metaphor, since no infinitesimals are actually used in such courses
except at best in some motivating discussions aimed to enhance student
intuitions.

In our course, two types of definitions of three key mathematical
concepts (continuity, uniform continuity, and convergence) were given:
\begin{enumerate}
\item[(A)] the usual~$\epsilon,\delta$ definition;
\item[(B)] the infinitesimal definition.
\end{enumerate}

The course first presented the infinitesimal definition (B-track, for
\emph{Bernoullian}) and then the~$\epsilon,\delta$ definition
(A-track, for \emph{Archimedean}).  We amplified the treatment in
Keisler due to the demands of the second semester sequel taught the
EDC way.  In particular, we expanded Keisler's treatment of the
$\epsilon,\delta$ approach, and added a treatment of the concept of
uniform continuity.  The following points should be kept in mind.
\begin{enumerate}
\item
The first edition of Keisler's book, \cite{Ke71}, was the first ever
calculus textbook using rigorous infinitesimals.
\item
Neither Archimedes nor Bernoulli envisioned anything like the
set-theoretic ontology underpinning the construction of modern
punctiform continua.  Modern terms like \emph{Archimedean continuum}
and \emph{Bernoullian continuum} refer not to \emph{ontology} of
mathematical entities but rather to the \emph{procedures} typically
used in the respective frameworks.
\item
The procedures in Robinson's framework provide closer proxies for the
procedures of historical infinitesimalists like Leibniz, Bernoulli,
Euler, and Cauchy than do the procedures in the modern Weierstrassian
framework, which similarly relies on a punctiform continuum.
\item
A popularisation of infinitesimals exploiting the field of rational
functions was developed by D. Tall under the name \emph{superreal
number system} (similarly a punctiform continuum); see e.g.,
\cite{Ta79}.  However this system lacks a transfer principle (see
Section~\ref{s23}) and cannot serve as a basis for a rigorous course
in the calculus.
\end{enumerate}
The $\epsilon,\delta$ definitions were a triumph of formalisation
mathematically speaking, but create pedagogical difficulties when
introduced without preparation, according to most scholars who have
studied the problem; see e.g., \cite{Da14}.  Our approach enables the
teacher to prepare the students for $\epsilon,\delta$ by explaining
the concepts first using a rigorous infinitesimal approach.  Studies
of methodology involving modern infinitesimals include \cite{Su76},
\cite{He04}, \cite{El}, \cite{PS}, \cite{TK}, \cite{Vi16}.

We sought to impart the fundamental concepts of the calculus in a way
that is the least painful to the students, while making sure that they
have the necessary background in the~$\epsilon,\delta$ techniques to
continue in the second semester course taught via EDC.  Once the
students understand the basic concepts via their intuitive B-track
formulations, they have an easier time relating to the A-track
paraphrases of the definitions.

Recently Robinson's framework has become more visible thanks to
high-profile advocates like Terry Tao; see e.g., his work \cite{Ta14},
\cite{TV}.  The field has also had its share of high-profile
detractors like Errett Bishop and Alain Connes.  Their critiques were
critically analyzed in \cite{KK11d}, \cite{KL13}, and \cite{KKM}.
Further criticisms by J. Earman, K. Easwaran, H. M. Edwards, Ferraro,
J. Gray, P. Halmos, H. Ishiguro, G. Schubring, and Y. Sergeyev were
dealt with respectively in the following recent texts: 
\cite{KS1}, \;
\cite{Ba14}, \; 
\cite{KKKS}, \; 
\cite{Ba17}, \, 
\cite{Bl16c}, \,
\cite{Bl16a}, 
\cite{Ba16b},
\cite{Bl16d}, 
\cite{Gu16}.
In \cite{BK} we analyze the Cauchy scholarship of Judith Grabiner.
For a fresh look at Simon Stevin see \cite{KK12}.

\section{Review of definitions}

In this section we will review both the A-track and the B-track
definitions of the three calculus concepts that the students were
polled on.  

\subsection{Procedure vs ontology}

As a prefatory remark, we would like to respond to a common objection
that to do calculus with infinitesimals you need to get a PhD in
nonstandard analysis first, and that's obviously not a good way of
teaching freshman calculus.

One way of responding is that a certain amount of foundational
material needs to be taken for granted in either approach.  Thus, the
real number system is not constructed in the EDC approach.  Instead,
certain subtle properties, like the existence of limit, closely
related to the completeness of the reals, are assumed on faith.

There is general agreement that in a calculus course we don't
elaborate the exact details concerning the real numbers with respect
to their ontological status in foundational theories such as the
Zermelo-Fraenkel set theory (ZFC).  This is because the
\emph{procedures} of the calculus don't depend on the
\emph{ontological} issues of set-theoretic axiomatisations.  In the
traditional approach to the calculus, we present all the procedures
rigorously, including all the epsilon-delta definitions, while staying
away from such ontological and foundational issues.

Similarly in developing TIC, we don't elaborate the set-theoretic
issues of the precise ontological status of the hyperreals in a ZFC
framework.  Rather, we teach them the procedures of the calculus
exploiting infinitesimals in a fully rigorous way, including the more
intuitive infinitesimal definitions of the key concepts.

\subsection{A-track definitions of key concepts}

A function~$f$ is said to be \emph{continuous} at a point~$c\in\R$ if
the following condition is satisfied:
\begin{equation}
\label{e21a}
(\forall\epsilon\in\R^+)(\exists\delta\in\R^+)(\forall x\in\R)\big[
|x-c|<\delta \;\Rightarrow\; |f(x)-f(c)|<\epsilon\big].
\end{equation}

A function~$f$ is said to be \emph{uniformly continuous} in a
domain~$D\subseteq \R$ if the following condition is satisfied:
\[
\begin{aligned}
(\forall\epsilon\in\R^+) & (\exists\delta\in\R^+) (\forall x\in D)
(\forall x'\in D) \\&
\big[|x'-x|<\delta\;\Rightarrow\;|f(x')-f(x)|<\epsilon\big].
\end{aligned}
\]

A sequence~$(u_n)$ is said to \emph{converge} to~$L\in\R$ if the
following condition is satisfied:
\[
(\forall\epsilon\in\R^+)(\exists N\in\N) (\forall n\in\N) \big[n>N
\;\Rightarrow\; |u_n-L|<\epsilon\big].
\]
The above A-track definitions are succinct summaries of the distilled
mathematical content of these concepts; of course, they were not
presented in this fashion in freshman calculus, but rather in a
slower-paced fashion.

\subsection{B-track definitions of key concepts}
\label{s23}

We now review the corresponding B-track definitions in more detail,
since they are less likely to be familiar to modern readers usually
educated in EDC frameworks.  What is involved is a hyperreal
extension~$\R\hookrightarrow\astr$, where~$\astr$ is an ordered field
including both infinitesimal (see below) and infinite numbers.  A key
tool in working with such an extension is the \emph{transfer
principle} (see below).

Such fields can be constructed from sets of sequences of real numbers,
similar to the construction of the real numbers from the rational
numbers.

An infinitesimal~$\alpha\in\astr$ is a number satisfying~$|\alpha|<r$
for every positive real~$r$.  An infinite number~$H$ is a number
satisfying~$|H|>r$ for every real number~$r$.

We teach the students to work with these new numbers, and to apply the
basic rules of arithmetic to them (infinitesimal times infinitesimal
is infinitesimal, infinitesimal times infinite can have any order of
magnitude, etc.), and what ``being infinitely close" means (see
Section~\ref{s24}).

We then introduce the \emph{standard part function} (sometimes called
\emph{shadow}).  This is a function from the finite (i.e., not
infinite) hyperreals to the reals, which rounds off each finite
hyperreal to its nearest real number.

The \emph{transfer principle} is a type of theorem that, depending on
the context, asserts that rules, laws or procedures valid for a
certain number system, still apply (i.e., are ``transfered'') to an
extended number system.  Thus, the familiar
extension~$\Q\hookrightarrow\R$ preserves the properties of an ordered
field.  To give a negative example, the
extension~$\R\hookrightarrow\R\cup\{\pm\infty\}$ of the real numbers
to the so-called \emph{extended reals} does not preserve the
properties of an ordered field.  The hyperreal extension
$\R\hookrightarrow\astr$ preserves \emph{all} first-order properties.
For example, the identity~$\sin^2 x+\cos^2 x=1$ remains valid for all
hyperreal~$x$, including infinitesimal and infinite values
of~$x\in\astr$.  Another example of a transferable statement is the
property that
\begin{quote}
for all positive reals~$x,y$, if~$x<y$ then~$1/y<1/x$.  
\end{quote}
Transfer applies to formulas like \eqref{e21a} that quantify over
elements of~$\R$, but not directly to statements that quantify over
\emph{sets} of elements.  Thus, the completeness property of the
reals, which involves quantification over sets, does not transfer
directly.  For a more detailed discussion, see the textbook
\emph{Elementary Calculus} \cite{Ke86}.

\subsection{Microcontinuity}
\label{s24}

Both continuity and uniform continuity can be defined in terms of an
auxiliary concept of \emph{microcontinuity} (this term is not used by
Keisler but it is used in \cite{Da77}).  The definition of
microcontinuity exploits the natural extensions~$\astf$ of a real
function~$f$ and~${}^\ast\!  D$ of a real set~$D$, available in a
hyperreal setting.

Let~$D=D_f\subseteq\R$ be the domain of a real function~$f$.  We say
that~$\astf$ is \emph{microcontinuous} at~$x$ if
\begin{equation}
\label{e21}
\text{whenever~$x'\approx x$, one also has~$\astf(x') \approx
\astf(x)$,}
\end{equation}
for all~$x'$ is in the domain~${}^\ast\!  D_{\!f}\subseteq\astr$
of~$\astf$.  Here the relation~$\approx$ is the relation of infinite
proximity, i.e., ~$x'\approx x$ if and only if the difference~$x'-x$
is infinitesimal.

The condition of microcontinuity can be tested not only at a real
point~$x\in D_f$ but also at an arbitrary hyperreal point~$x\in {}
^\ast\! D_{\!f}$.

Thus, the squaring function
\begin{equation}
\label{e23}
y=x^2
\end{equation}
fails to be microcontinuous at an infinite point~$H\in\astr$.  Indeed,
let~$\alpha=\frac{1}{H}$ so that $H\approx H+\alpha$.  Note that the
corresponding~$y$-increment is
\[
\Delta y= (H+\alpha)^2-H^2=H^2+2H\alpha+\alpha^2-H^2= 2+\alpha^2.
\]
It follows that the~$y$-increment is appreciable rather than
infinitesimal, showing that the squaring function is not
microcontinuous at~$H$.

It turns out that a real function~$f$ is \emph{continuous} at~$c\in\R$
if and only if~$\astf$ is microcontinuous at~$c$, and
\emph{uniformly continuous} in its domain~$D_f$ if and only if the
following condition is satisfied:
\[
(\forall x\in{}^\ast\!D_f)(\forall x'\in{}^\ast\!D_f) \left[ x\approx
x' \; \Rightarrow \; \astf(x)\approx \astf(x') \right].
\]
Equivalently, a real function~$f$ is uniformly continuous on~$D_f$ if
$\astf$ is microcontinuous at~$x$ for each~$x \in {}^\ast\!D_f$.

To continue with the example \eqref{e23} presented above, one can now
state that the failure of the squaring function to be uniformly
continuous on~$\R$ is due to its failure to be microcontinuous at a
single infinite point.  This proof of the failure of uniform
continuity of a function is of reduced quantifier complexity when
compared to A-track proofs of the same fact.

\subsection{Convergence}

The last of the three concepts we tested is \emph{convergence}.  A
sequence~$(u_n)$ converges to~$L$ if and only if
\[
\text{st}(u_H)=L
\]
for all infinite values~$H$ of the index, where ``st'' denotes the
\emph{standard part function}.

Consider for example the sequence~$u_n=\frac{n+1}{n}$.  To prove that
the limit is~$1$, we write 
\[
\lim_{n\to\infty} u_n = \text{st}(u_H) =
\text{st}\left(\frac{H+1}{H}\right)=\text{st}\left(1+\frac{1}{H}\right)=1
\]
by the additive property of the shadow, since~$\frac{1}{H}$ is
infinitesimal.  Note that rather than having to deal with an
\emph{inverse} problem as in the EDC framework, the proof is a
\emph{direct} calculation.  Moreover, it is just as rigorous as the
EDC proof, the difference being that the EDC proof would necessarily
involve preliminary calculations or at least a guess for a limit
value.

\section{The poll}

The questionnaire, which mostly followed a multiple-choice format,
also contained a control question asking students to prove that
\[
\lim_{x\to 2}(x+5)=7
\]
in two different ways:~$\epsilon,\delta$ (A-track) and infinitesimal
(B-track); see Section~\ref{one}.

Almost all the students (98\%) attempted to solve the control problem
via B-track, while 71\% attempted to give an A-track solution.  Of
those who attempted a B-track solution, 85\% succeeded; of those who
attempted an A-track solution, 20\% succeeded.

   It should be noted that our students had substantial practice
specifically in~$\epsilon,\delta$.  Our TAs spent two entire sessions
on this, and the students also had to submit homework assignments
where they were required to use the~$\epsilon,\delta$ technique.
Also, they did many exercises similar to the above using
$\epsilon,\delta$.

The students were asked to comment on the helpfulness of A-track and
B-track definitions of three key concepts: continuity, uniform
continuity, and convergence of a sequence.  More specifically, they
were presented with the statement ``the definition helped me
understand the concept,'' and were given the following five options
for a possible answer: (1) agree strongly; (2) agree; (3) undecided;
(4) disagree; (5) disagree strongly.

With respect to the B-track definition of continuity, 69\% of the
students felt that the definition helped them understand the concept
(``agree'' or ``agree strongly'').  Meanwhile, 10\% of the students
felt that the A-track definition of continuity helped them understand
the concept.  Among students who were able to define continuity
correctly, 75\% felt the B-track definition helped them understand the
concept, while 9\% felt the A-track definition helped them understand
it.

With regard to uniform continuity, 74\% felt that the B-track
definition helped them understand the concept, whereas 21\% felt that
the A-track definition helped them understand the concept.  Among
students who were able to define uniform continuity correctly, 80\%
felt that the B-track definition helped them understand it, whereas
24\% felt that the A-track definition helped them understand it.

With regard to the definition of a convergent sequence, 62\% felt that
the B-track definition helped them understand the concept, whereas
10\% felt that the A-track definition helped them understand it.
Among students who were able to define convergence correctly, 70\%
felt that the B-track definition helped them understand it, whereas
13\% felt that the A-track definition helped them understand it.

\section{Divide-and-conquer vs paraphrase}
\label{s4}

In our poll, the percentage of students who felt that the B-track
definition helped them understand the concept increases by about 7\%
(of the respondents) when one calculates the percentage on the basis
of those students who were able to give a correct definition of the
appropriate concept.  A similar phenomenon occurs among students who
felt that the A-track definition helped them understand the concept in
the case of the concepts of uniform continuity and convergent
sequences.

On average among the three concepts, over two-thirds (68\%) of the
students felt that the B-track definition is helpful, while only about
one in seven students (14\%) felt that the A-track definition is
helpful.

To summarize, what we tried to do in the course is to impart to the
students the fundamental concepts of the calculus in a way that is the
least painful to the students, while making sure that they have the
necessary background in the~$\epsilon,\delta$ techniques to continue
in the second semester course taught via EDC.  The results of the poll
suggest that starting with the intuitive B-track definitions succeeds
in this sense.  Once the students understand the basic concepts via
their intuitive B-track formulations, they are able to relate more
easily to the A-track paraphrases of the definitions.

\subsection{From Pythagorians to Weierstrassians and beyond}

To comment on the idea of paraphrase in more detail, suppose one is
limited to working with the rational numbers (for fear of getting
thrown overboard by enraged Pythagorians).

Yet in one's mathematical investigations there may arise a need to
express the predicate that an unknown rational number~$x$ is greater
than the diagonal of a unit square.  However, one is only allowed to
use inequalities of the form
\[
x>q
\]
where~$q$ is rational.  Since one is forbidden to talk about
irrationals, one says instead that~$x$ is greater than every
rational~$q$ such that~$q^2 <2$, or in formulas
\begin{equation}
\label{e41}
\forall q\in\Q \left[ q^2<2 \implies q < x \right].
\end{equation}
This quantified formula looks more complicated than the intended
inequality, but since one already knows what it means, one can readily
understand it.  In other words, the complicated quantified
formula~\eqref{e41} is merely a long-winded paraphrase for the
familiar inequality $x>\sqrt{2}$.

Similarly, someone interested in property~\eqref{e21} that an
infinitesimal change in input should always produce an infinitesimal
change in output, may be led to exploit the $\epsilon,\delta$ formula
\eqref{e21a} with its notorious alternating quantifiers (to avoid
Hippasus' fate at the hands of enraged Weierstrassians).

These quantified formulas look complicated, but they are merely
long-winded paraphrases for simpler definitions exploiting
infinitesimals that were used by Cauchy but have been suppressed since
1870 when Weierstrass and his followers broke with the infinitesimal
mathematics of Leibniz, Euler, and Cauchy.

Our approach could therefore be described as an application of the
divide-and-conquer algorithm.  One first separates the inherent
difficulty of the subject of the calculus into two parts:
\begin{enumerate}
\item
the intrinsic difficulty of the concepts themselves;
\item
the technical complications of the A-track paraphrases with their
notorious quantifier alternations.
\end{enumerate}

The good news is that the concepts are accessible without the A-track
paraphrases thanks to \cite{He48}, \cite{Lo55}, and \cite{Ro61}.  The
idea is to start with (1), contrary to the EDC approach that starts
with (2).  Our approach is more consistent with Toeplitz's thinking,
discussed in the next section.

\section{All the stops out}

Otto Toeplitz had the following to say in 1927 about teaching
infinitesimal calculus:
\begin{quote}
I consider it an inviolable axiom that by the end of a two-semester
course, a beginner should have obtained a complete understanding and
complete mastery of the technique of `epsilontic' operations, and that
he did not bring such techniques with him from high school. \; \;
\cite[p.~303]{To27}
\end{quote}
We heartily agree with the latter point, and empathize with the
former.  Toeplitz continues:
\begin{quote}
The way this has been formulated already suggests the solution.
Instead of launching the `epsilontic' methodology right away at the
beginning with all the stops out, as one says of the organ, one should
lead the student gradually up a gentle ascent to the peak of this
technique, just as the organist uses one register after the other in a
well-composed piece of organ music - and in this way not one of the 45
percent spoken of above will be left out.  (ibid.)
\end{quote}
The explanation of Toeplitz's 45\% figure is as follows.  Toeplitz
considers that about~5\% of the students are the ``natural"
mathematicians that will grasp the \emph{Epsilontik} immediately and
do not even need to go to the lectures.

Toeplitz also considers that about 50\% of the student body present in
the mathematics courses is too weak, making it difficult to structure
the course based on them.  Dismissing half the class in this fashion
is unacceptable, and arguably is a consequence of an obligatory
adherence to the EDC approach.  The TIC approach makes it possible to
reach close to 100\% of the students, by postponing the introduction
of $\epsilon,\delta$ definitions as we explained in Section~\ref{s4}.

Toeplitz is talking about the remaining 45\% who are strong students.
He is arguing that the top 5\% should not be taught at the expense of
the 45\%.

Note that Toeplitz advocated using infinitesimals and differentials at
a time (1927) when it was still considered that they are irremediably
lost in a hazy fog of meaninglessness, as Courant colorfully put it in
his textbook.%
\footnote{Courant described infinitesimals on page 81 of Differential
and Integral Calculus, Vol I, as ``devoid of any clear meaning" and
``naive befogging.''  Similarly on page 101, Courant described them as
``incompatible with the clarity of ideas demanded in mathematics,''
``entirely meaningless,'' ``fog which hung round the foundations,''
and a ``hazy idea.''  Cantor, Russell, and the mathematicians of
Courant's generation were convinced that infinitesimals are
self-contradictory.  Following \cite{Ro61} we know this not to be the
case.}
Thus, in discussing Kepler's Second Law, Toeplitz does not hesitate to
exploit both (infinitesimal) differentials and the notion of
\emph{utter smallness} as a pedagogical device:
\begin{quote}
In Figure 128,~$dF$ is the area of a narrow sector of area formed by
two closely neighboring radii; it can be approximated by the right
triangle~$SPQ$ ($PQ$ perpendicular to~$r$ at~$P$). The little chip by
which it exceeds~$dF$ is \emph{utterly small} in relation to~$dF$,
and, as~$dF$ itself becomes smaller, this chip diminishes even more
rapidly and therefore can be neglected. \cite[p.~151]{To63} [emphasis
added]
\end{quote}
Toeplitz was hardly the only one to exploit the explanatory power of
such terms.  A majority of mathematics educators involved in teaching
calculus routinely exploit such expressions.  They would say things
like ``the function~$f$ has limit~$L$ at a point~$a$ if we can make
$f(x)$ as close to~$L$ as we wish for all~$x$ sufficiently close to
$a$.''  They would say ``given epsilon positive and as small as we
wish\ldots" These expressions are variations on Toeplitz's ``utterly
small.''

This is how a majority of educators explain things to students, and
this is the language they use, because this is the way we think and
the way we perceive these notions.  It is not merely a pedagogical
device, but this is how we understand these ideas.

TIC offers us a possibility of making our intuitions precise with
terms like ``infinitesimal" (in place of ``as small as we wish") and
``infinitely close" (instead of ``as close as we want").

Many mathematicians think in terms of Toeplitz's ``utter smallness''
and related ideas.  The TIC approach makes effective use not only of
the students' intuitions but also of the mathematicians' intuitions
about infinitesimals.

\section*{Acknowledgments}

We are grateful to our students for their stimulating participation in
both the course and the poll.  We thank Jacques Bair and Piotr B\l
aszczyk for helpful comments on an earlier draft of the paper.
M.\;Katz was partially supported by the Israel Science Foundation
grant no.\;1517/12.

\end{document}